\documentclass[12pt]{amsart}
\input amssym.def
\input amssym.tex
\usepackage{amsfonts, amsmath, amsthm}
\usepackage{psfig}

\newtheorem{pro}{Proposition}[section]
\newtheorem{Theo}[pro]{Theorem}
\newtheorem{Lemma}[pro]{Lemma}

\newtheorem{cor}[pro]{Corollary}

\theoremstyle{definition}

\theoremstyle{remark}

\begin{document}

\begin{abstract}
  We define the \textit{injectivity radius} of a Coxeter polyhedron in
 $\textbf{H}^3$ to be half the shortest translation length among
 hyperbolic/loxodromic elements in the orientation-preserving reflection
 group.  We show that, for finite-volume polyhedra, this number is always
 less than 2.6339... , and for compact polyhedra it is always less
 than 2.1225... .
\end{abstract}

\title[Injectivity radii of hyperbolic polyhedra]
{Injectivity radii of hyperbolic polyhedra}

  \author{Joseph D. Masters}
 \maketitle


\section{Introduction}

Hyperbolic reflection groups are discrete groups of isometries of hyperbolic
space generated by reflections in the faces of a polyhedron.
  They provided some of the earliest known examples of Kleinian groups, and
 have been well-studied (see [V]).  In this paper, we prove that
 3-dimensional hyperbolic reflection groups always contain short
 elements. 

To be more precise, let $\Gamma^+(P)$ be the group of orientation-preserving
 isometries generated by reflections in the faces of the polyhedron $P$.  We
 define
 $injrad(\Gamma^+(P))$ to be half the shortest translation length
 among hyperbolic/loxodromic elements of $\Gamma^+(P)$.  Then we have: 

\vspace{1pc}

\textbf{Theorem 4.1}
\textit{(Main Theorem)  Let $P$ be a finite-volume Coxeter polyhedron in} $\textbf{H}^3$.
\textit{Then $injrad(\Gamma^+(P)) < cosh^-1(7)$ = 2.6639... . If $P$ is compact, then} \break
\textit{$injrad(\Gamma^+(P)) <  cosh^-1(3+4cos(2\pi/5))$ = 2.1225...} .

\vspace{1pc}

Remarks:\\
1. It is known (see [W]) that if $\lbrace\textbf{H}^3/ \Gamma_i \rbrace$ is
 a family of closed hyperbolic 3-manifolds,
 and if $\lbrace rank(\Gamma_i) \rbrace$ is bounded, then
 $\lbrace injrad(\textbf{H}^3/ \Gamma_i) \rbrace$ is also bounded
  (recall the \textit{rank} of a finitely generated group is the cardinality
 of a minimal generating set).  Observe that by [B], $Rank(\Gamma(P))$
 increases with the number of sides of $P$, so Theorem 4.1
 is not covered by [W].\\
\\
2.  For more about short geodesics in hyperbolic
 3-manifolds, see [AR].\\
\\
3. We speculate that the bounds may be sharp, but we do not know a proof.\\
\\

\textbf{Idea of Proof}:
  To prove the Main Theorem, we must show that every three-dimensional
 hyperbolic reflection group contains a hyperbolic element with suitably
 short translation length.  This can usually be done by finding two
 non-adjacent faces 
 of the polyhedron which are suitably close; the short element is obtained
 by composing the reflections in the corresponding hyperplanes.

   A result of Nikulin's guarantees that every Coxeter polygon which is not
 a triangle has two non-adjacent sides which are close, and a two-dimensional
 version of the Main Theorem follows easily.  
 For most polyhedra, we can use Nikulin's result to show that
 two non-adjacent faces are close; the exceptions are those which contain
 ``non-prismatic'' faces (see Section 2 for a definition).  We show that
 these exceptional cases always contain triangular faces.  Then, after
 extending to the sphere at infinity, we use combinatorics and Euclidean
 geometry to deduce the existence of a short element.  The results
 on hyperbolic polyhedra contained in Section 3 allow us to sharpen the bound.

\textbf{Organization}
  Section 2 contains basic definitions; Section 3 contains some general,
 technical results about hyperbolic polyhedra; Section 4 contains the
 proof of the Main Theorem.

\vspace{1pc}

\textbf{Acknowledgments}:
  I would like to thank Dr. Alan Reid for his help and encouragement.
  Thanks also to Dr. Gary Hamrick and David Bachman for valuable conversations,
 Dr. Daryl Cooper for helpful correspondence, and the referee for
 valuable comments.

\section{Definitions}

A \textit{convex polyhedron}, $P$,  in $\textbf{H}^n$ is a countable intersection of
 closed half-spaces: $P = \bigcap_i H_i^+$, where $H_i$ denotes a hyperplane
 and $H_i^+$ the corresponding closed half-space.
When n=2, we use the term $polygon$ instead.
If $P$ is a convex polyhedron in $\textbf{H}^n$, we let $\Gamma(P)$ denote the group of isometries generated by reflections 
in the bounding hyperplanes of $P$.  $\Gamma^+(P)$ denotes its orientation-preserving subgroup of index 2.
We say that a finite-volume, convex  polyhedron $P$ is a \it{Coxeter polyhedron }\rm  
if its dihedral angles are all integer submultiples of $\pi$;
if $P$ is a Coxeter polyhedron, then $\Gamma^+(P)$ is discrete.  
Two faces of $P$ are \textit{adjacent} if they share an edge.
  Given a hyperplane $H$ in $\textbf{H}^n$,
 $\rho_H$ will denote the isometry obtained by reflection in $H$.
  We will denote the hyperbolic distance between two sets $X$ and $Y$
 in $\textbf{H}^n$ by $d(X,Y)$.

  Given an n-sided face $F$ of $P$ $(n > 3)$, label the 
 edges of $F$ by $E_1, ..., E_n$, where $E_i$ shares a vertex
 with $E_{i+1}$ for $i = 1, 2, ..., n-1$,
 and label the
 adjacent faces by $F_1, ... , F_n$, where $F$ and $F_i$ share edge $E_i$
 (we say the faces adjacent to $F$ are labeled ``cyclicly'').
  We say that $F$ is $\textit{prismatic}$ if, for $i,j = 1,...,n$,  $|i-j| > 1
 \, (mod\, n)$ implies $F_i$ is non-adjacent to $F_j$.  

  By a \textit{hyperbolic n-manifold}, $\textbf{H}^n/ \Gamma$, we mean the
 quotient of hyperbolic n-space by a discrete group of
 isometries acting freely.  If $\Gamma$ has torsion, the quotient
 space $H^n/\Gamma$ is a \textit{hyperbolic n-orbifold}.
 The \it{injectivity radius}  \rm  of a hyperbolic manifold
 $M = \textbf{H}^n/ \Gamma$ is equal to $sup \lbrace \alpha \in \mathbb{R}^+ |$
 \textit{every point x} $\in$ \textit{M is the center of an embedded ball of
 radius} $\alpha\rbrace$.
  We shall generalize this
 definition to the case where $\Gamma$ has torsion:
 the \textit{injectivity radius of a Kleinian group} $\Gamma$, denoted
 $injrad(\Gamma)$, is equal
 to half the shortest translation length among hyperbolic/loxodromic elements
 of $\Gamma$.  Note that this agrees with the usual notion when $\Gamma$
 is torsion-free.  The \textit{injectivity radius of a Coxeter polyhedron} $P$
 is equal to $injrad(\Gamma^+(P))$.
 Given an element $g$ of
 $\Gamma$, we will denote its translation length by $\ell(g)$.

\section{Hyperbolic Polyhedra}

The following technical result shall be used in our proof of the Main Theorem.

\begin{Theo}
Let P be a finite-volume, convex polyhedron in $\textbf{H}^3$ with acute dihedral angles
 and no triangular faces.  Then P has a prismatic face.  If, furthermore, P
 is compact, then P has a prismatic quadrilateral or pentagonal face.
\end{Theo}
\begin{proof}
We first consider the finite-volume case.

\begin{Lemma}
  Let P be a finite-volume, convex polyhedron in
 $\textbf{H}^3$ with acute dihedral angles and no triangular faces, and suppose that P contains at least
 one non-prismatic face.  
 Let $G$ be a planar graph representing the 1-skeleton of $P$ (see Fig. 1).
 Then there are three non-prismatic
 faces of P which bound a region in G consisting entirely of prismatic faces.
\end{Lemma}

\begin{proof}
Let $F$ be a non-prismatic face of P, and label its adjacent faces cyclicly
 by $F_1,F_2,...$.  We have that $F_i$ and $F_j$ are adjacent for some i,j with
$|i-j| > 1$ (see Fig. 1).
   Note that $F$, $F_i$, and $F_j$ bound a region $R$,
 and that $F_i$ and $F_j$ are also non-prismatic.
  Since $P$ contains no triangular faces, $R$ cannot be a face.
  And if the sublemma is false, $R$ must contain a non-prismatic face $F'$.
  Label the faces adjacent to $F'$ cyclicly by $F_1', F_2',...$.
  Then for some $k,\ell$ with $|k-\ell| > 1$, $F_k'$ and $F_{\ell}'$
 are adjacent.
  Then $F'$, $F_k'$, and $F_{\ell}'$
 bound a triangular region $R' \subsetneq R$.
  Again, if the sublemma is false, $R'$ must contain a non-prismatic
 face $F''$,
 creating a triangular region $R'' \subsetneq R'$. 
  Since $P$ has a finite number of sides, this process must eventually
 terminate with three non-prismatic faces bounding a region containing
 only prismatic faces.  This proves the sublemma, and the non-compact
 case of Lemma 3.1 follows immediately.

\end{proof}

\begin{figure}[htbp]

\begin{center}

\ \psfig{file=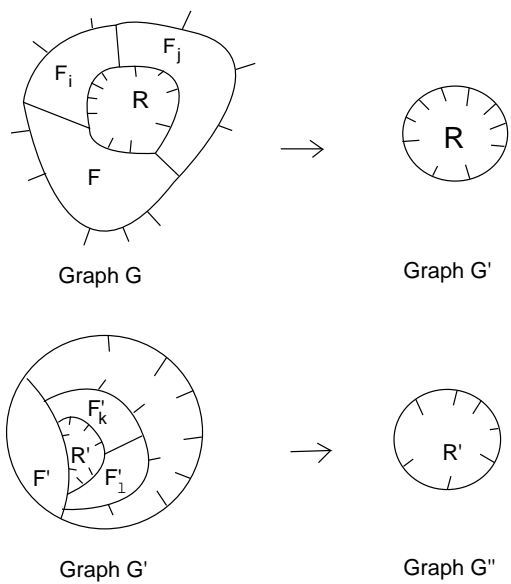,height=3in}

\caption{Reducing the graph around a non-prismatic face}

\end{center}

\end{figure}

To prove the stronger statement for the compact case, we shall need the
 following combinatorial lemma:

\begin{Lemma}
Let $P$ be a compact, convex polyhedron in $\textbf{H}^3$ with acute dihedral angles. Then $P$ contains a face
 with 5 or fewer edges.
\end{Lemma}

\begin{proof}
  This follows from an Euler characteristic count.
  Let $|V|$ = number of vertices of $P$, $|E|$ = number of edges of $P$,
 and $|F|$ = number of faces of $P$.  Since $P$ is a compact, convex hyperbolic
 polyhedron, $P$ is \textit{simple}-- i.e. each vertex is shared by
 exactly 3 different edges (see [A1]).  So  $|V| = 2|E|/3$.

  We have  $|V| - |E| + |F| = 2$.

 $|F| - |E|/3 = 2$.

  $|F|(1 - |E|/3|F|) = 2$.

  So $|E|/|F| < 3$.   

So the average number of edges per face $<$ 6.  So $P$ must contain a face with
 5 or fewer edges.

\end{proof}

Now suppose $P$ is compact.  If all faces of $P$ are prismatic, we are done
 by Lemma 3.3.  So suppose $P$ contains a non-prismatic face.  Then by Lemma 3.2,
 there are three non-prismatic faces of $P$ which bound a region $R$ in $G$
 consisting entirely of prismatic faces.  We need to show that $R$ contains a
 face with 5 or fewer edges.  We form $\hat{R}$ from $R$ by subtracting the
 three vertices, $v_1,v_2,v_3$, of non-prismatic faces on the boundary of $R$,
 and then adding a face to the exterior of $R$, so that
 $\hat{R}$ is topologically a sphere (see Fig. 2).

\begin{figure}[htbp]

\begin{center}

\ \psfig{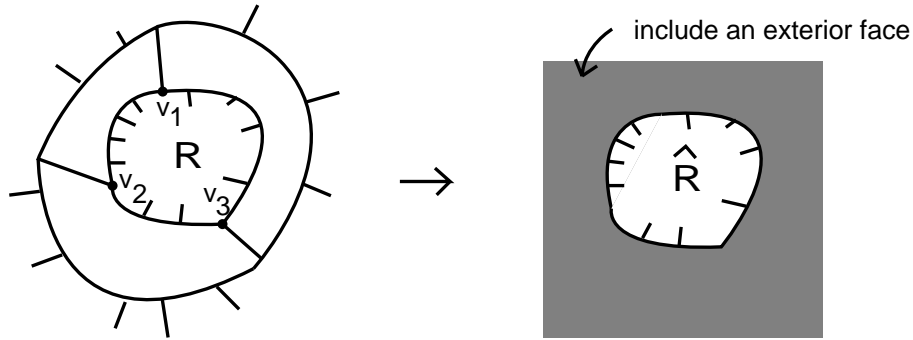}

\caption{To form $\hat{R}$ from $R$, we remove three vertices and include an exterior face}

\end{center}

\end{figure}

  Then, arguing as in Lemma 3.3,
 $\hat{R}$ must
 contain a face $F$ with 5 or fewer edges.  We claim that $R$ must also contain
 such a face.  This will certainly be true unless

\vspace{1pc}

1.  $F$ is the face on the exterior of $R$.   

\vspace{.5pc}
or
\vspace{.5pc}

2.  $F$ is a pentagon, and $F$ contains one of the edges from which a vertex
 was deleted (note that $F$ can contain at most one of the $v_i$'s, since it
 is prismatic).
\vspace{1pc}

If $F$ is the face on the exterior and $F$ is a triangle, then we claim that
 $R$ must consist of three quadrilaterals, as in Fig. 3b.  For otherwise,
 $R$ must contain three non-prismatic faces, contrary to assumption (see Fig. 3b).   
 
\begin{figure}[htbp]

\begin{center}

\ \psfig{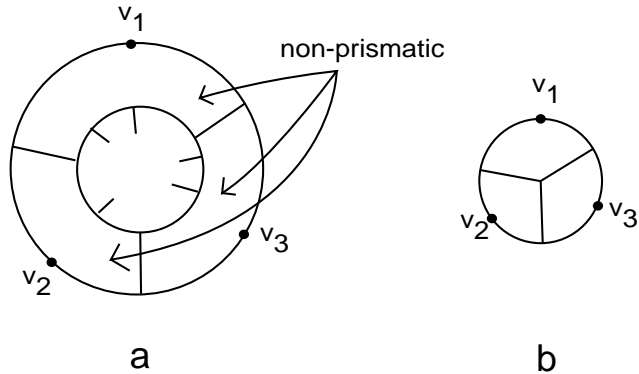}

\caption{If the exterior is a triangle, then $R$ will contain either
 a) non-prismatic faces, or b) three quadrilaterals.}

\end{center}

\end{figure}

So if the theorem is false, then $\hat{R}$ contains at most four faces
 with fewer than 6 sides; at least three of them must be pentagons,
 and none are triangles.  We shall show that this is impossible.

Number the faces of $\hat{R}$ by $F_1, F_2, ..., F_n$. Let $|E_j|$ be the number
 of sides of $F_j$.  Then $|E_1| + ... + |E_n| = 2|F|$, where $|F|$ is
 the total number of faces of $\hat{R}$.  Since $P$ is simple, the graph is
 trivalent, so $|F|$ is divisible by 3; hence $|E_1| + ... +|E_n|$ is
 divisible by 6.
  Since also the average number of sides per face is $<$ 6, we have:
    
$6n-5 \leq \sum_{i=1,..,n}|E_i| < 6n$.

Thus the sum is not divisible by 6, for a contradiction.

\end{proof}

\section{Injectivity radius of hyperbolic polyhedra} 

In this section we prove the Main Theorem:
\begin{Theo}
Let P be a finite-volume Coxeter polyhedron in $\textbf{H}^3$.  Then \break
$injrad(\Gamma^+(P)) < acosh(7) = 2.6339...$ .  If P is compact, then \break
 $injrad(\Gamma^+(P)) < acosh(3 + 4cos(2\pi/5)) = 2.1226...$.
\end{Theo}

An important part of the proof of Theorem 4.1 is played by the following
 2-dimensional result.

\begin{Theo}
Let P be a finite-area Coxeter n-gon in $\textbf{H}^2$, $n > 3$.  Then \break
 $injrad(\Gamma^+(P)) \le acosh(3 + 4cos(2\pi/n)) < acosh(7) = 2.6339... $.  If n=3, then \break
 $injrad(\Gamma^+(P)) \le acosh(3 + 4cos(2\pi/4))$.  If P is compact,
 then all the inequalities are strict.
\end{Theo}
\begin{proof}
The theorem is a consequence of the following lemma, which is a re-phrasing of
 ([N], Theorem 3.2.1).

\begin{Lemma}
Let P be a finite area, convex n-gon in $\textbf{H}^2$.  Label the
bounding geodesics of P cyclicly by $H_1, H_2, ..., H_n$.  Then
for some i, $d(H_i, H_{i+2}) \le acosh(3 + 4cos(2\pi/n))$ (subscripts taken
 mod n).  If P is compact, the inequality is strict.
\end{Lemma}

\begin{proof}
We will reproduce Nikulin's proof that  $d(H_i, H_{i+2}) < acosh(7)$.
  The proof of the finer estimate is a bit more complicated, and we omit it
 (see [N]).

We use the Lobachevsky model for $\textbf{H}^2$.  First, pick a point $p$ on the interior of $P$.
Project the vertices of $P$ onto the circle at infinity along rays emanating from $p$.
This determines an ideal polygon $P'$ with bounding geodesics $H_1', H_2', ...$.
We claim that $d(H_i, H_{i+2}) \le d(H_i', H_{i+2}')$.  For the distance between
$H_i'$ and $H_{i+2}'$ is measured along a mutually orthogonal geodesic segment
 $\alpha$ (see Fig. 4a); since $\alpha$ must intersect $H_i$ and $H_{i+2}$, 
$d(H_i, H_{i+2}) < length(\alpha) = d(H_i', H_{i+2}')$. 

\begin{figure}[htbp]

\begin{center}

\ \psfig{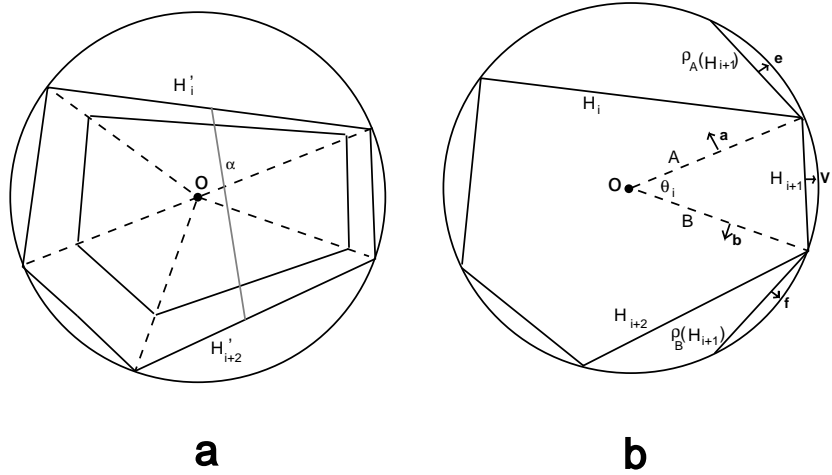}

\caption{a) The faces of the ideal polygon are farther apart.  b) Notation for Lemma 3.2}

\end{center}

\end{figure}

So it is enough to consider the case where $P$ is an ideal polygon.
Pick $i$ such that the Euclidean angle $\theta_i = v_iOv_{i+1}$ is minimal,
 where $O$ is the origin.
Let $A$ be the diameter through $v_i$ and let $B$ be the diameter through
$v_{i+1}$.  Then $d(H_i, H_{i+2}) \le d(\rho_A(H_{i+1}), \rho_B(H_{i+1}))
 = acosh(-<\bf{e}, \bf{f}>) $, where $\bf{e}$ is the unit normal to $\rho_A(H_{i+1})$, $\bf{f}$ is the unit normal to $\rho_B(H_{i+1})$,
and $<.|.>$ is the inner product
 $<(x_1,y_1,z_1)|(x_2,y_2,z_2)> = x_1x_2 + y_1y_2 - z_1z_2$ (see Fig. 4b).

Now, let $\bf{v},\bf{a}$ and $\bf{b}$ be outward unit normals to $H_{i+1}$,
 $A$ and $B$, respectively.  Then we have:

$<\bf{e},\bf{f}> = <\bf{v} - 2<\bf{v}, \bf{a}>\bf{a}, \bf{v} - 2<\bf{v}, \bf{b}>\bf{b}>$

\vspace{1 pc}

$= <\bf{v}, \bf{v}> - 2<\bf{v}, \bf{a}>^2 - 2<\bf{v}, \bf{b}>^2 + 4<\bf{v},\bf{a}> <\bf{v}, \bf{b}> <\bf{a}, \bf{b}>$ 

\vspace{1 pc}

$= 1 - 2 - 2 - 4cos(\theta_i)$.

\vspace{1 pc}

So $d(H_i,H_{i+2}) < acosh(7)$.
\end{proof}

We now resume the proof of Theorem 4.2.
Let $P$ be a Coxeter n-gon, and suppose first $n > 3$.
Pick two non-adjacent edges such that the
corresponding geodesics $H$ and $H'$ are less than $\break$
 $acosh(3 + 4cos(2\pi/n))$ apart.  Let
 $g = \rho_H\rho_{H'}$ in $\Gamma^+(P)$.  Since $P$ has acute angles, $H$ is
 disjoint from $H'$, so $g$ is hyperbolic.  Let $\alpha$ denote the axis of
 $g$, and note that it is perpendicular to both $H_i$ and $H_{i+2}$
 (see Fig. 5). $\ell(g)$ is given by $d(p, g(p))$, where $p$ is any point
 on $\alpha$. Taking $p$ to be $\alpha \cap H_{i+2}$, it is easy to see that
 $\ell(g) = 2 d(H_i, H_{i+2})$.
  Therefore $injrad(\Gamma^+(P)) \le \ell(g)/2 = d(H,H') <
 acosh(3 + 4cos(2\pi/n))$. 

\begin{figure}[htbp]

\begin{center}

\ \psfig{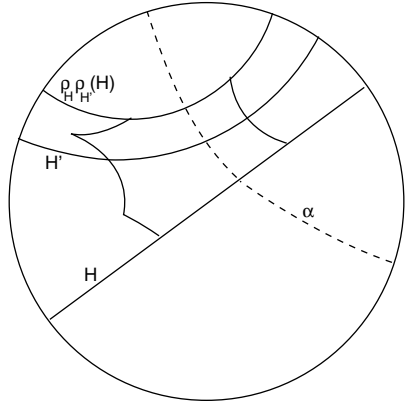}

\caption{The distance from $H$ to $H\prime$ is half the translation legth of $\rho_H \rho_{H\prime}$}

\end{center}

\end{figure}

If $P$ is a triangle, suppose first that $P$ has no right angles, and is
 compact.  Label the vertices
of $P$ by $v_0, v_1, v_2$, with corresponding angles
 $\theta_0, \theta_1, \theta_2$ (see Fig. 6a). 
One of the angles, say $\theta_0$, must be $\le \pi/4$.
Consider the quadrilateral
$Q$, with angles $(\theta_0, 2\theta_1, \theta_0, 2\theta_2)$,
 obtained by reflecting $P$ along $v_1v_2$.  Let $H$ and $H'$ be nonadjacent bounding geodesics
of $Q$, with $d(H,H') < acosh(3 + 4cos(2\pi/4))$.  $H$ and $H'$ are disjoint:  for if they intersected,
they would create a triangle $T$ with angles $\break
 (\pi - \theta_0, \pi - 2\theta_i, x)$, for
 i = 1 or 2; however $\theta_0 + 2\theta_i < \pi$, so $T$ would have angle sum $> \pi$, which
is impossible. 
Therefore, as above, 
$g = \rho_H\rho_{H'}$ is a hyperbolic element of $\Gamma^+(P)$ with $\ell(g)/2 < acosh(3 + 4cos(2\pi/4))$. 

\begin{figure}[htbp]

\begin{center}

\ \psfig{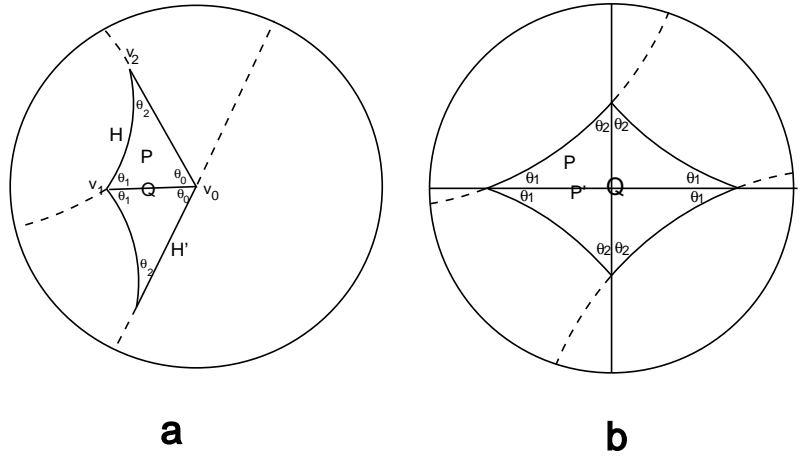}

\caption{We can reflect to obtain a quadrilateral.}

\end{center}

\end{figure}

If $P$ has a right angle at a vertex $v_0$, and is compact,
reflect P along $v_0v_1$ to obtain a new triangle, $P'$, with angles
 $(2\theta_1, \theta_2, \theta_2)$ (see Fig. 6b).
Reflect $P'$ along $v_0v_2$ to obtain a quadrilateral $Q$ with angles
 $(2\theta_1, 2\theta_2, 2\theta_1, 2\theta_2)$.  Then, again, the opposite geodesics
 bounding $Q$ must be disjoint, or else they would create a triangle with angles
 $(\pi - 2\theta_1, \pi - 2\theta_2, x)$, which is impossible since
 $\theta_1 + \theta_2 < \pi/2$.  So again we obtain a hyperbolic element $g$ with
 $\ell(g)/2 < acosh(3 + 4cos(2\pi/4))$.

If $P$ has an ideal vertex $v_0$, then reflecting $P$ along $v_1v_2$ creates a quadrilateral $Q$ 
whose bounding geodesics clearly cannot intersect, so that again we obtain
 the required hyperbolic element.
\end{proof}

As a corollary of Theorem 4.2, we have:

\begin{cor}
Let P be a finite-volume Coxeter polyhedron in $\textbf{H}^3$ with an n-sided prismatic face.
  Then $injrad(\Gamma^+(P)) \le acosh(3 + 4cos(2\pi/n))$.
\end{cor}

\begin{proof}
Consider the n-sided prismatic face $F_0$, and label the faces adjacent to
 $F_0$ cyclicly by $F_1, ..., F_n$.  Let $H_j$ denote the hyperplane spanned by
 $F_j$, and let $\break$
 $H_{i,j} = H_i \cap H_j$.
By Lemma 4.3, there is some i such that $\break$
$d(H_{0,i}, H_{0,i+2}) \le acosh(3 + 4cos(2\pi/n))$.
Then $d(H_i, H_{i+2}) \le acosh(3 + 4cos(2\pi/n))$.
  Since $F_i$ and $F_{i+2}$ are nonadjacent, $H_i$ and $H_{i+2}$ do
not intersect by [A1].  So $g = \rho_{H_i}\rho_{H_{i+2}}$ is a hyperbolic
 element of $\Gamma^+(P)$.   So
 $injrad(\Gamma^+(P)) \le \ell(g)/2 \le acosh(3 + 4cos(2\pi/n))$.   
\end{proof}

We now prove the Main Theorem.

\vspace{1pc}

\textit{Proof of 4.1.}

\vspace{1 pc}
\textbf{Compact case}:

\vspace{1 pc}
Case 1: \textit{$P$ has no triangular faces}.

\vspace{1 pc}
By Theorem 3.1, $P$ contains a prismatic quadrilateral or pentagonal face,
and so by Corollary 4.6, we are done.

\vspace{1 pc}

Case 2:  \textit{$P$ has a triangular face, but $P$ is not a simplex}.

\vspace{1 pc}
 Let $F_0$ be the triangular face, and label the faces adjacent to $F_0$
 by $F_1, F_2$ and $F_3$.  
   Let $P\prime = \bigcap_{i=0,1,2,3}H_i^+$, where $H_i$ denotes the
 hyperplane spanned by $F_i$ (recall Section 2 for the definition of $H_i^+$).
  Let $C_i$ and $C_i^+$
 denote the boundary at infinity of $H_i$ and $H_i^+$, respectively.
  Let $\theta_{ij} = $ dihedral angle between $F_i$ and $F_j = \pi/n_{ij}$ 
 (here we are again using the fact that in the compact case
 the polyhedra are simple).  
 We label the dihedral angles of $P\prime$ by
  $\textbf{A}_{P\prime}\rm = ((n_{01},n_{02},n_{03}),(n_{12}, n_{23}, n_{31}))$. 
 Then, 

\vspace{1 pc}

 I. \hspace{1in}       $1/n_{0i} + 1/n_{0j} + 1/n_{ij} > 1$, and 

\vspace{1 pc}

II.  \hspace{1in}     $1/n_{12} + 1/n_{23} + 1/n_{31} < 1$.

\vspace{1 pc}

One may then easily verify that the only 5 possibilities for
 $(n_{01},n_{02},n_{03})$ are: (2,2,2), (2,2,3), (2,2,4), (2,2,5) and (2,3,3).

\vspace{1 pc}

Case 2a:  $\textbf{A}_{P\prime} = ((2,2,2),(n_{12}, n_{23}, n_{31}))$
 or $((2,2,4),(n_{12}, n_{23}, n_{31}))$.

\vspace{1pc}
Let $Stab_{F_0}$ denote the subgroup of $\Gamma^+(P)$ which leaves $F_0$
invariant.  Then Theorem 5.4 of [BM] implies that $Stab_{F_0}$ contains
 a triangle group.  So by Theorem 4.2, we are done.

\vspace{1 pc}
Case 2b:   $\textbf{A}_{P\prime} = ((2,3,3),(n_{12}, n_{23}, n_{31}))$.

\vspace{1pc}
By I and II, the only two possibilities (modulo relabeling of edges) are $\break$
 $\textbf{A}_{P\prime} = ((2,3,3),(4,2,5))$ or $((2,3,3),(5,2,5))$.
  Consider $P\prime\prime = P\prime \cup \rho_{H_1}(P\prime)$ (see Fig. 7a).
  Let $Q$ denote the quadrilateral created in $H_0$.  Conjugate so that,
 in the upper half space model, $C_1$ is the imaginary axis (see Fig. 7b),
 and $\infty
 \in C_i^+-C_i$ for i = 2,3.  $H_1\cap H_2\cap H_3 = \emptyset$, since
 the angles $\theta_{0i}$ are acute;  therefore, $\break$
   $C_2^+ \cap C_3^+ \cap \rho_{H_1}(C_2^+) \cap \rho_{H_1}(C_3^+)
 = Q_\infty^1 \cup Q_\infty^2$,
 where $Q_\infty^1, Q_\infty^2$ are quadrilateral regions in
 $S^2_{\infty}= \hat{\mathbb{C}}$,
and $\infty \in Q_\infty^2$ .  The dihedral angles of $Q_\infty^1$
 are $(\theta_{23}, 2\theta_{12}, \theta_{23}, 2\theta_{31})$.
 
\begin{figure}[htbp]

\begin{center}

\ \psfig{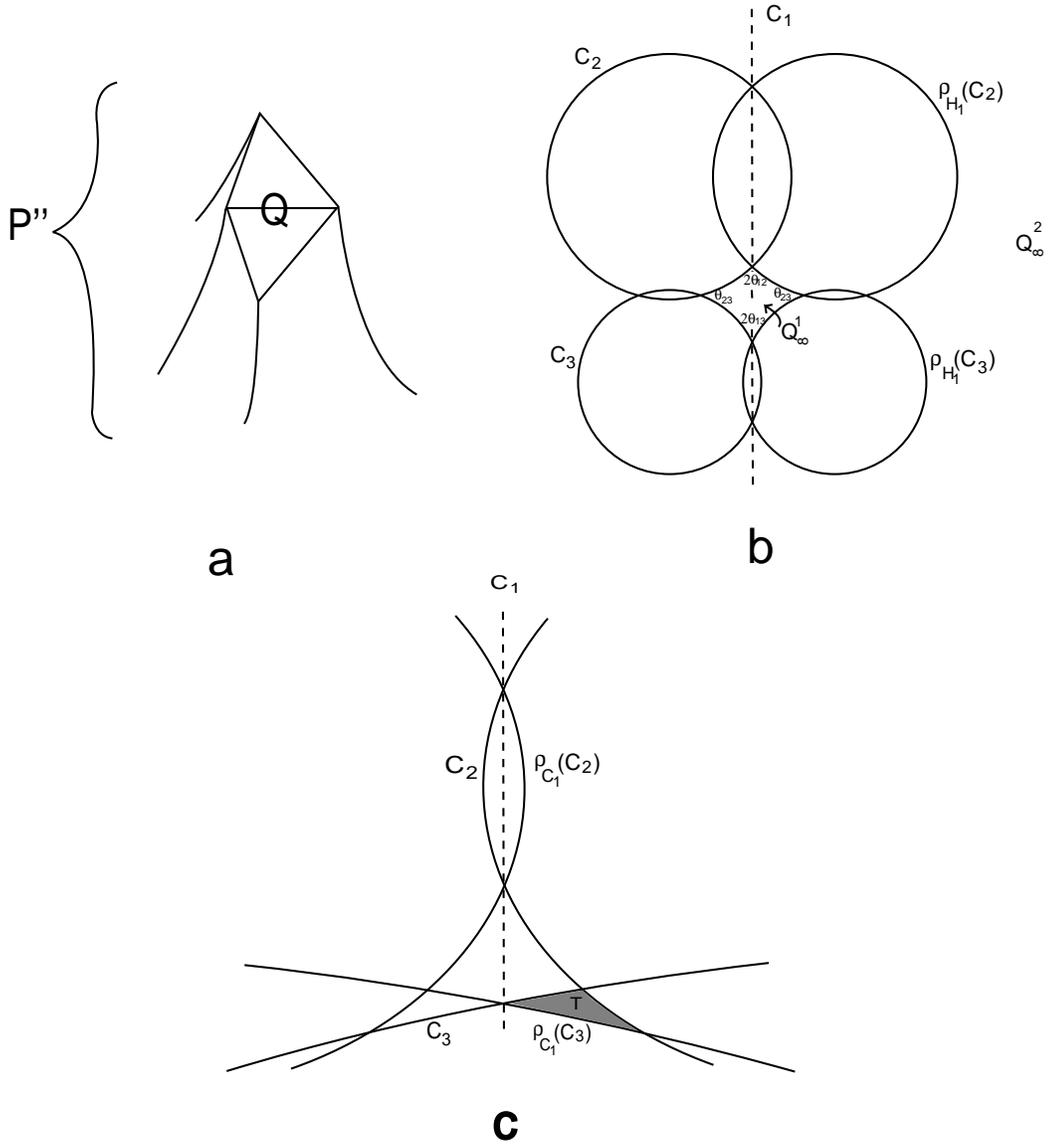}

\caption{a) The polyhedron $P\prime\prime$  b) The view on the sphere at infinity c) If opposite circles intersect, a triangle is formed.}
\end{center}

\end{figure}

 We claim that $C_2 \cap \rho_{H_1}(C_3) = C_3 \cap \rho_{H_1}(C_2)
 = \emptyset$, in which case $H_2, \rho_{H_1}(H_3)$ are disjoint, and
 $H_3, \rho_{H_1}(H_2)$ are disjoint.
  Otherwise, three circles would intersect to create a
 triangle $T$ with angle sum $ > \pi$ (see Fig. 7c).
  However, $T$ has only one positively curved side ---
 call its curvature $\kappa$ --- and it has another, longer side of curvature
 $-\kappa$.  So, by the Gauss-Bonnet formula, $T$ has angle sum $< \pi$, for
 a contradiction.  Hence
 $g_1 = \rho_{H_2}\rho_{\rho_{H_1}(H_3)}$
 and $g_2 = \rho_{H_3}\rho_{\rho_{H_1}(H_2)}$ are both hyperbolic elements of
 $\Gamma^+(P)$.  By Lemma 4.3, two of the opposite faces of $Q$ must be
 less than $acosh(3 + 4cos(2\pi/4))$ apart.  It follows that either $g_1$ or
 $g_2$ has suitably short translation length.

\vspace{1 pc}

Case 2c:  $\textbf{A}_{P\prime} = ((2,2,3),(n_{12}, n_{23}, n_{31}))$ or
  $((2,2,5),(n_{12}, n_{23}, n_{31}))$.

\vspace{1 pc}
Without loss of generality, assume $n_{13} \ge n_{23}$, so by II,
 $n_{13} \ge 3$.
  First suppose $n_{12}, n_{23} \ge 3$.  
  As in case 2b, reflect in $H_1$ to create a polyhedron
 $P\prime\prime$ with a quadrilateral face $Q$.  Again, on $S_{\infty}^2$ we see
 a quadrilateral $Q_\infty^1$ with angles
 $(\theta_{23}, 2\theta_{12}, \theta_{23}, 2\theta_{31})$.  
  Since the sum of any two adjacent angles of $Q_\infty^1$ is $\le \pi$,
 we can argue as in case 2b to show that opposite circles bounding $Q\prime$
 must be disjoint, thus creating a
 hyperbolic element with suitably short translation length.

  If $n_{12} = 2$, then by I and II, $n_{23}$ and $n_{31} \ge 4$.  Then after
 reflecting in $H_1$ and in $H_2$, we see on $S_{\infty}^2$ an
 acute quadrilateral (see Fig. 8), and we may argue as above.

\begin{figure}[htbp]

\begin{center}

\ \psfig{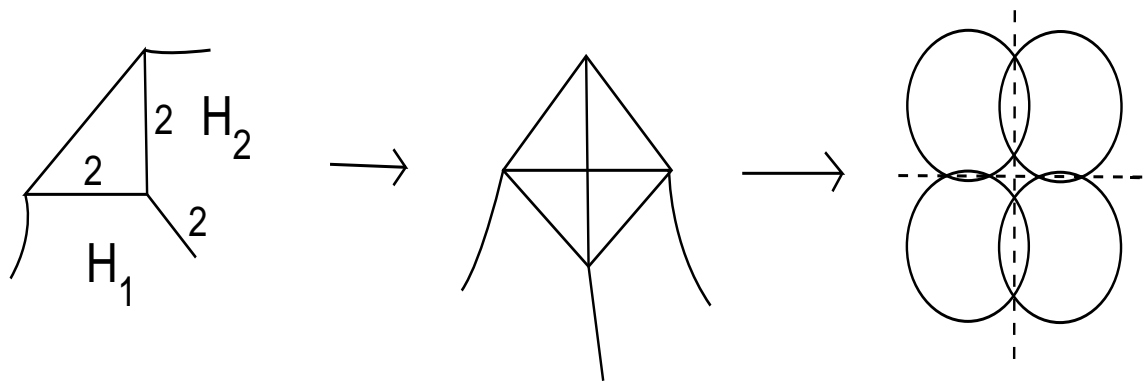}

\caption{Reflect twice, then view at infinity.}

\end{center}

\end{figure}

  If $n_{23} = 2$, then by I and II, $n_{12} \ge 4$.  
    
  ---if $n_{12} = 4$, then
 after reflecting in $H_2$, we can reduce to the case where $n_{12} =2$.

  ---if $n_{12} \ge 5$, then by II, $n_{23} \ge 4$, and reflecting in $H_2$
 creates an acute quadrilateral.  So we may argue as in case 2b.

\vspace{1 pc}
Case 3: \textit{P is a simplex}.

\vspace{1 pc}
By [L] there are only nine congruence classes of compact simplices in $\textbf{H}^3$.
  By [M] (see also [BM]), eight of these contain triangle groups,
 so in these cases the result follows
 from Theorem 4.3.  Denote the remaining tetrahedron by $T_8$;
 label its faces $F_1, .., F_4$; and let $\pi/n_{ij}$ be the dihedral angle
 between $F_i$ and $F_j$.  We have $n_{12}=2,\, n_{13}=3,\, n_{14}=4,\, n_{23}=5,
 \, n_{24}=3$, and $n_{34}=4$.  
  Let $H_i$ be the hyperplane spanned by $F_i$.  It is not
 difficult to construct $T_8$ explicitly in Lobachevsky space
 and then compute the faithful discrete representation of $\Gamma^+(T_8)$ in
 $O(3,1)$.  Then it is straightforward to compute that
 $\rho_{H_3}\rho_{H_4}\rho_{H_2}\rho_{H_1}\rho_{H_4}\rho_{H_2}$ is a
 hyperbolic element with translation length
 1.66131... $< 2acosh(3 + 4cos(\pi/5))$.

\vspace{2 pc}

\textbf{Non-compact case:}

\vspace{1 pc}

By Theorem 3.1 and Corollary 4.6, it is enough to consider the case where
 $P$ has a triangular face, $F_0$.  As in the compact case, the faces
 adjacent to $F_0$ form a polyhedron $P\prime$.

  If none of the vertices
 of $F_0$ are ideal, then the proof for the compact case carries over without
 change. So suppose $F_0$ has an ideal vertex.

 Ideal vertices may be either tri-valent or 4-valent (see [A2]).
  If the vertices of $F_0$ are all trivalent then condition II still holds;
 condition I holds at regular vertices, and changes to an equality at
 ideal vertices.
 Then the techniques from the compact case
 are sufficient to prove the theorem-- we omit the details.
  If one of the vertices is 4-valent,
 then two of the sides adjacent to $F_0$ are tangent on $S^2_{\infty}$
 (see Fig. 9).
  One may now argue as in the compact case, treating the tangent sides as
 adjacent with dihedral angle 0.  The theorem follows.

\begin{figure}[htbp]

\begin{center}

\ \psfig{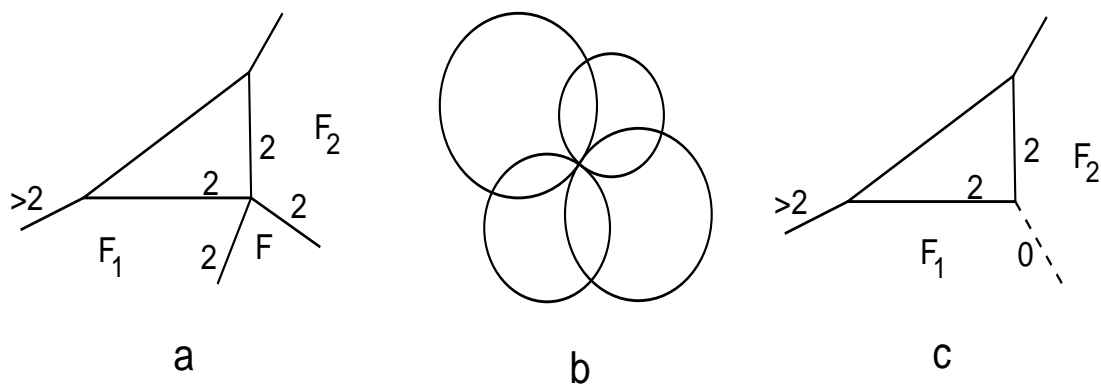}

\caption{a) A triangular region with a 4-valent ideal vertex.  b) A 4-valent
 ideal vertex on $S_{\infty}^2$.  c) We may remove $F$ and view $F_1$ and
 $F_2$ as adjacent with dihedral angle 0.}

\end{center}

\end{figure}

$\qed$

\section{References}

[A1] E. M. Andreev, \textit{On convex polyhedra in
 Loba$\check{c}$evski$\breve{i}$ spaces}, Math. USSR-Sb. \textbf{10} (1970),
 413-440.

\vspace{1pc}

[A2] E. M. Andreev, \textit{On convex polyhedra of finite volume in
 Loba$\check{c}$evski$\breve{i}$ spaces}, Math. USSR-Sb. \textbf{12} (1970),
 No. 2, 255-259. 

\vspace{1pc}

[AR] C. C. Adams and A. W. Reid, ``Systoles for hyperbolic knot and
 link complements'',  Math. Proc. Camb. Phil. Soc., to appear.

\vspace{1pc}

[B]  K. Brown, ``Buildings'', Springer-Verlag, New York, 1989.

\vspace{1pc}

[BM]  T. Baskan and A.M. Macbeath, \textit{Centralizers of reflections in crystallographic groups}, Math. Proc. Camb. Phil. Soc. \textbf{92} (1982), 79-91.

\vspace{1pc}

[L]  F. Lanner, \textit{On complexes with transitive groups of automorphisms},
 Comm. Sem. Math. Univ. Lund 11 (1950).

\vspace{1pc}

[M]  C. Maclachlan, \textit{Triangle subgroups of hyperbolic tetrahedral
 groups}, Pacific J. Math. \textbf{176} (1996) 195-203.

\vspace{1pc}

[N]  V. V. Nikulin, \textit{On arithmetic groups generated by reflections in
 Loba$\check{c}$evski$\breve{i}$ spaces}, Math. USSR Izv. \textbf{16} (1981), 573-601.

\vspace{1pc}

[V]  E. B. Vinberg, \textit{Hyperbolic reflection groups}, Russian Math. Surveys \textbf{40} (1985), 31-75.

\vspace{1pc}

[W]  M. White, pre-print.

\vspace{2pc}

Department of Mathematics

University of Texas at Austin 78712

masters@math.utexas.edu

\vspace{2pc}

Current Address

Department of Mathematics

MS 136, Rice University

6100 S. Main St.

Houston TX 77005-1892

\end{document}